\documentclass[11pt]{amsart}
\usepackage{amsmath,amssymb,amsthm}
\usepackage{hyperref}
\hypersetup{hidelinks=true}

\usepackage{bbm}

\usepackage[utf8]{inputenc}

\usepackage{graphicx}
\usepackage{enumerate}
\usepackage{amsfonts,amssymb,latexsym}
\usepackage{accents}
\usepackage{afterpage}
\usepackage[style=alphabetic,maxbibnames=99, backend=bibtex]{biblatex}
\usepackage{tikz}

\newcommand{\N}{\mathbb{N}}

\newcommand{\R}{\mathbb{R}}

\newcommand{\s}{\sum}

\newcommand\restr[2]{{
  \left.\kern-\nulldelimiterspace 
  #1 
  \vphantom{\big|} 
  \right|_{#2} 
  }}

\newtheorem{theorem}{Theorem}[section]
\newtheorem{lemma}[theorem]{Lemma}

\newtheorem{prop}[theorem]{Proposition}
\newtheorem{corollary}[theorem]{Corollary}
\theoremstyle{definition}
\newtheorem{definition}[theorem]{Definition}

\theoremstyle{remark}
\newtheorem{remark}[theorem]{Remark}
\numberwithin{equation}{section}

\def\fnote#1{\footnote}

\def\R{{\mathbb R}}

\def\ignora#1{}
\def\n3#1{\left\vert  \! \left\vert \! \left\vert \, #1 \, \right\vert \!
  \right\vert \! \right\vert }

\include{draft}

\begin{document}



\title[Bases in nets]{Schauder bases in Lipschitz free spaces over nets of $\mathcal{L}_{\infty}$-spaces}

\author{Petr H\'ajek}\thanks{This research was supported by CAAS CZ.02.1.01/0.0/0.0/16-019/0000778
 and by the project  SGS21/056/OHK3/1T/13.}
\address[P. H\'ajek]{Czech Technical University in Prague, Faculty of Electrical Engineering.
Department of Mathematics, Technická 2, 166 27 Praha 6 (Czech Republic)}
\email{ hajek@math.cas.cz}

\author{ Rub\'en Medina}\thanks{The second author research has also been supported by MICINN (Spain) Project PGC2018-093794-B-I00 and MIU (Spain) FPU19/04085 Grant.}
\address[R. Medina]{Universidad de Granada, Facultad de Ciencias.
Departamento de An\'{a}lisis Matem\'{a}tico, 18071-Granada
(Spain); and Czech technical University in Prague, Faculty of Electrical Engineering.
Department of Mathematics, Technická 2, 166 27 Praha 6 (Czech Republic)}
\email{rubenmedina@ugr.es}
\urladdr{\url{https://www.ugr.es/personal/ae3750ed9865e58ab7ad9e11e37f72f4}}

\maketitle

\begin{abstract}
In the present note we give a construction (based on a retractional argument) of a Schauder basis for the Lipschitz free space $\mathcal{F}(N)$, 
over a net $N$ in any  separable infinite dimensional $\mathcal{L}_\infty$-space $X$. In particular, this yields  the first example of an infinite dimensional Banach space $X$ not containing $c_0$ with such a property.
\end{abstract}

\section{Introduction}

The topic of our note was motivated by the important and well-known question of Kalton \cite{Kal12} whether every (separable) Banach space is approximable, or equivalently whether the Lipschitz free space $\mathcal{F}(N)$ over a net $N$
in any separable Banach space has the Bounded Approximation Property (BAP for short).
Kalton proved that for a Banach space $X$ with a separable dual (or a separable dual space itself), as well as for Banach spaces $X$ with the BAP,  the answer is positive.  

In the present note we study the retractional structure of nets in Banach spaces $X$ of the type $\mathcal{L}_\infty$.
As a corollary of our construction we obtain the result that the Lipschitz free space $\mathcal{F}(N)$ over a net
in any separable $\mathcal{L}_\infty$-space has a Schauder basis. This is a much stronger condition than the BAP,
guaranteed by Kalton's results, and moreover our proof is entirely constructive. It should be noted that our main result
extends our previous papers \cite{HM21}, resp. \cite{HN17}, where a retractional structure of nets in Banach spaces with a Schauder
basis (resp. FDD) containing a copy of $c_0$ has been investigated (for the same purpose).
One of the key ingredients in our previous results was the coefficient quantization 
theorem of \cite{DOS+08}, which requires that the Banach space $X$ contains a copy of $c_0$.
 Our present 
results are a genuine extension, as there exist (by a famous result of Bourgain and Delbaen \cite{BD80})  separable $\mathcal{L}_\infty$-spaces without a copy of $c_0$. Although the retractional approach
to the problem is the same as in our previous papers, the construction is rather different
and relies instead on the explicit description of the abstract $\mathcal{L}_\infty$-spaces as
Bourgain-Delbaen spaces obtained in \cite{BD80}, and an alternate quantization which crucially uses the ambient $\ell_\infty$-space.
It should be added that the class of $\mathcal{L}_\infty$-spaces is very significant and in some sense large, and extremely useful
for constructing (counter)examples (\cite{BD80}, \cite{AH11}). In particular,  every Banach space with a separable dual is
contained in an isomorphic $\ell_1$-predual (of the type $\mathcal{L}_\infty$) \cite{FOS11}.

Let us now proceed with some technical background.

Given $a,b\in\R^+$ we will say that a subset $N$ of a metric space $M$ is an $(a,b)$-net of $M$ whenever it is $a$-separated and $b$-dense. More precisely, $N$ is considered to be $a$-separated whenever $d(x,y)\ge a$ for every $x,y\in N$, $x\neq y$, and $b$-dense if for every $x\in M$ there is $y\in N$ such that $d(x,y)\le b$. We will say that $N\subset M$ is a net whenever there are $a,b\in \R^+$ such that $N$ is an $(a,b)$-net.

A map $T$ from a metric space $M$ into another metric space $N$ is said to be Lipschitz if there exists some $\lambda>0$ such that
$$d(T(x),T(y))\le \lambda d(x,y) \;\;\;\;\forall x,y\in M.$$
We say that $\lambda$ is a Lipschitz constant for $T$ and we call the infimum of all Lipschitz constants for $T$ the Lipschitz norm of $T$, that is,
$$||T||_{\text{Lip}}=\sup\limits_{x,y\in M, x\ne y}\frac{d\big(T(x),T(y)\big)}{d(x,y)}.$$
If $\lambda>0$ is a Lipschitz constant for $T$ then we say that $T$ is $\lambda$-Lipschitz.

\begin{definition}
Let $N_1$ and $N_2$ be metric spaces. We say that $N_1$ and $N_2$ are Lipschitz equivalent whenever there is a bijection $F:N_1\to N_2$ such that both $F$ and $F^{-1}$ are Lipschitz.
We put $\text{dist}(F)=\|F\|_{\text{Lip}}\|F^{-1}\|_{\text{Lip}}$ the Lipschitz distorsion of $F$.
\end{definition}

In this note, we are mainly interested in finding a retractional structure in nets. To this end, we are going to order the elements of the net
into a sequence and introduce $K$-Lipschitz retractions onto the first $n$ points of the sequence with a commuting behaviour. More precisely,


\begin{definition}
Let $M$ be a countable metric space and $K\ge1$, a $K$-retractional basis of $M$ is a sequence of retractions $\{\varphi_n\}_{n=1}^\infty$, $\varphi_n:M\to M$ satisfying,
\begin{enumerate}
\item $\varphi_n(M)=M_n:=\bigcup_{j=1}^{n}\{x_j\}$,
\item $\bigcup_{j=1}^\infty\{x_j\}=M$,
\item $\varphi_n$ is $K$-Lipschitz for every $n\in\N$,
\item $\varphi_m\circ\varphi_n=\varphi_{\min(m,n)}$ for every $m,n\in\N$.
\end{enumerate}
We will say that $M$ has a retractional basis whenever it has a $K$-retractional basis for some $K\ge1$.
\end{definition}

It is clear that these retractions can be used to extend Lipschitz maps from a finite number of points $\{x_1,\dots,x_n\}$ to the whole net by composition, keeping the Lipschitz constant under control independently of the number of points from which the extension is made.

The existence of retractional bases has been studied previously (see \cite{HN17}, \cite{HM21} and \cite{Nov20}).
It is easy to see that retractional bases in nets are clearly preserved under Lipschitz equivalences. More precisely,

\begin{prop}\label{lipequiv}
Let $N_1$ and $N_2$ be countable Lipschitz equivalent metric spaces, with distorsion $D$. If $N_1$ has a $K$-retractional basis for some $K>0$ then $N_2$ has a $DK$-retractional basis. 
\end{prop}

A Schauder basis for a real Banach space $X$ is a sequence $(x_n)\subset X$ with the property that for every $x\in X$, there exists a unique sequence $(\alpha_n)\subset \R$ such that
$$\bigg|\bigg| x-\sum\limits_{i=1}^n\alpha_ix_i \bigg|\bigg|\xrightarrow{n\to\infty}0.$$
If $(x_n)$ is a Schauder basis then there is a $K>0$ such that $\Big|\Big|\sum\limits_{i=1}^n\alpha_ix_i \Big|\Big|\le K||x||$ for every $x\in X$ and $n\in\N$. In this case we say that $(x_n)$ is a $K$-Schauder basis.

A retractional basis for a net implies that the Lipschitz-free space over the net has a Schauder basis, made out of the linearization of the retractions.

\begin{theorem}[\cite{HN17}]\label{theoretr}
Let $M$ be a countable metric space with a $K$-retractional basis $\{\varphi_n\}_{n=1}^\infty$. Then the Lipschitz-free space over $M$ has a $K$-Schauder basis.
\end{theorem}

In \cite{HM21} we obtained a retractional basis for nets in Banach spaces with a Schauder basis and containing a copy of $c_0$, resp. having a $c_0$-like FDD. It is known \cite{JRZ71} that $\mathcal{L}_\infty$-spaces have a Schauder basis, but as shown in \cite{BD80} they do not necessarily contain a copy of $c_0$.

For background on the approximation properties of Banach spaces we refer to the Handbook article \cite{Cas01}. We will use the standard notation and terminology of the Banach spaces theory, as in \cite{Fab1}.
For background on Lipschitz-free spaces and its main properties we refer to \cite{GK03}.

\section{Nets in $\mathcal{L}_\infty$-spaces}

The $\mathcal{L}_\infty$-spaces were introduced by Lindenstrauss and Pelczynski in \cite{LP68}.
They have played an important role in Banach space theory, and their structure theory
is   a very actively pursued topic of research (see for example \cite{LR69},\cite{JRZ71},\cite{BD80},\cite{BP83},\cite{AH11} and \cite{AGM16}). A crucial step in this theory was provided by Bourgain and Delbaen \cite{BD80} 
 who introduced an inductive construction scheme of such spaces with various prescribed properties.
In particular, it was shown that $\mathcal{L}_\infty$-spaces need not necessarily contain a copy of $c_0$.
It turned out recently \cite{AGM16}, that in fact the scheme can be applied to all separable $\mathcal{L}_\infty$-spaces. We will rely on this approach in our paper.

\begin{definition}
A Banach space $X$ is called an $\mathcal{L}_{\infty,\lambda}$-space, for some $\lambda\ge1$, if for every finite dimensional subspace $F$ of $X$ there is a finite dimensional subspace $G$ of $X$, containing $F$ and $\lambda$-isomorphic to $\ell_\infty^n$ for $n=\text{dim}\,G$. $X$ is said to be an $\mathcal{L}_\infty$-space if it is an $\mathcal{L}_{\infty,\lambda}$-space for some $\lambda\ge1$.
\end{definition}

Throughout this section, we are going to follow the notation, terminology and concepts used in \cite{AGM16}. In particular, we are going to use the equivalent definition of $\mathcal{L}_\infty$-spaces given in \cite{AGM16}.

\begin{definition}
Let $\Gamma_1,\Gamma$ be non-empty sets with $\Gamma_1\subset \Gamma$. A linear operator $i: \ell_\infty(\Gamma_1)\to\ell_\infty(\Gamma)$ is said to be an extension operator if for every $x\in \ell_\infty(\Gamma_1)$ and $\gamma\in \Gamma_1$ we have that $x(\gamma)=i(x)(\gamma)$.
\end{definition}

\begin{definition}
Let $(\Gamma_n)_{n\in\N}$ be a strictly increasing sequence of non-empty sets and $\Gamma=\bigcup_{n\in\N}\Gamma_n$. A sequence of extension operators $(i_n)_{n\in\N}$ with $i_n:\ell_\infty(\Gamma_n)\to\ell_\infty(\Gamma)$ for every $n\in\N$ are said to be compatible if for every $n,m\in\N$ with $n<m$  it holds that $i_n=i_m\circ r_m\circ i_n$, where $r_n:\ell_\infty(\Gamma)\to\ell_\infty(\Gamma_n)$ denotes the natural restriction operator.
\end{definition}

\begin{definition}\label{BDspace}
A Banach space $X$ is said to be a Bourgain-Delbaen space if there is a strictly increasing sequence $(\Gamma_n)_{n\in\N}$ of non-empty finite sets and a sequence of compatible extension operators $(i_n)_{n\in\N}$ with $i_n:\ell_\infty(\Gamma_n)\to \ell_\infty(\Gamma)$ where $\Gamma=\bigcup_{n\in\N}\Gamma_n$ such that $\lambda=\sup_n||i_n||<\infty$ and 
$$X=\overline{\bigcup\limits_{n\in\N}i_n(\ell_\infty(\Gamma_n))}.$$
\end{definition}

It is proved in Theorem 3.6 of \cite{AGM16} that every separable infinite dimensional $\mathcal{L}_\infty$-space is isomorphic to a Bourgain-Delbaen space. The reverse result is clearly satisfied. From now on, let $X$ be a Bourgain-Delbaen space. Let us denote by $(\Gamma_n)$ and $(i_n)$ the sequences given in definition \ref{BDspace}, and $\lambda=\sup_n||i_n||$. We will also denote $\Gamma=\bigcup_{n\in\N}\Gamma_n$.

For every $r\in\R$ we define the entire part of $r$ as 
$$[r]=\frac{r}{|r|}\max\big\{s\in\N\cup\{0\}\;:\;s\le|r|\big\}$$
if $r\neq 0$ and $[0]=0$.
Given a subset $S$ of $\Gamma$, we define the quantization of the elements of $\ell_\infty(S)$  by $f:\ell_\infty(S)\to\ell_\infty(S)$ where
$$f(x)(\gamma)=[x(\gamma)],$$
for every $x\in \ell_\infty(S)$ and $\gamma\in S$. Also, we define for every $s>0$ the truncation of radius $s>0$ in $\ell_\infty(S)$ as $T_s:\ell_\infty(S)\to sB_{\ell_\infty(S)}$ given by
$$T_s(x)(\gamma)=\begin{cases}x(\gamma)\;\;&\text{if }|x(\gamma)|\le s,\\s\frac{x(\gamma)}{|x(\gamma)|}&\text{if }|x(\gamma)|>s.\end{cases}$$
If $\Gamma_n\subset S$ for $n\in \N$ we will denote by $r_n:\ell_\infty(S)\to\ell_\infty(\Gamma_n)$ the natural restriction operator.
It is worth mentioning that $r_n$, $f$ and $T_s$ are nonexpansive. In dealing with $f$, $r_n$ and $T_s$ we simplify the notation by not mentioning the set $S$, since it will be clear from the context.

Our intention is to prove that every net in $X$ has a retractional basis. The proof is  splitted  into six steps. In the first two steps we define a separated set $M\subset \ell_{\infty}(\Gamma)$ and retractions onto some finite subsets of $M$. In the steps 3 to 5 we construct the retractional basis for $M$. Finally, in the last step we prove that there is a net in $X$ that is Lipschitz equivalent to $M$. This finishes the proof since any two nets of an infinite-dimensional Banach space are necessarily Lipschitz equivalent \cite{LMP00}.

\subsection{Step 1}

In the first step we are going to define the set $M$ by describing an increasing sequence of finite sets $(M_n)_{n\in\N}$ whose union will be $M$. We will also construct a sequence of commuting Lipschitz retractions $(\phi_n)$ from $M$ onto $M_n$.

First, for every $n\in\N$ we set $\widetilde{M}_n=i_n\circ f\big( \lambda^nB_{\ell_{\infty}(\Gamma_n)} \big)$. Then, we define inductively an increasing sequence of sets $(M_n)\subset \ell_\infty(\Gamma)$ by taking $M_1=f\big(\widetilde M_1\big)$ and for every $n\in\N$,
$$M_n=M_{n-1}\cup f\circ i_n\Big(f\big( \lambda^nB_{\ell_{\infty}(\Gamma_n)} \big)\setminus r_n\big( M_{n-1} \big)\Big).$$

Now we define $M$ as $\bigcup\limits_{n\in\N}M_n$, which is clearly $1$-separated. In order to construct the Lipschitz retractions onto the blocks $M_n$ we need the following Lemmas.

\begin{lemma}\label{dense}
For every $n\in\N$ and $x\in M_n$ the following inequality holds
$$\big|\big| i_n\circ r_n(x)-x \big|\big|\le \lambda+1.$$
\begin{proof}
We let $\widetilde x=i_n\circ r_n(x)\in \widetilde M_n$. Now, let us take $m\le n$ such that $x\in M_m\setminus M_{m-1}$ (considering $M_0=\emptyset$). Then, by the definition of $M_m\setminus M_{m-1}$ and the compatibility condition on the extensions it follows that
$$x=f\circ i_m\circ r_m(x)=f\circ i_m\circ r_m(\widetilde x)=f\circ i_n\circ r_n\circ i_m\circ r_m(\widetilde x).$$
Then, let
$$z:=f(\widetilde x)=f\circ i_n\circ r_n(x)=f\circ i_n\circ r_n\circ f\circ i_m\circ r_m(\widetilde x).$$
Finally,
$$\begin{aligned}||\widetilde x-x||=&||\widetilde x-z||+||z-x||\\=&||\widetilde x-f(\widetilde x)||\\&+||f\circ i_n\circ r_n\circ f\circ i_m\circ r_m(\widetilde x)-  f\circ i_n\circ r_n\circ  i_m\circ r_m(\widetilde x)\big)||\\\le&1+||f\circ i_n\circ r_n||_{\text{Lip}}||f\circ i_m\circ r_m(\widetilde x)-i_m\circ r_m(\widetilde x)||\le1+\lambda.\end{aligned}$$
\end{proof}
\end{lemma}

\begin{lemma}\label{Lipmain1}
Given $n\in\N$, the map $\restr{r_n}{M_n}$ is an injection. Moreover, $r_n(M_n)=f(\lambda^n B_{\ell_\infty(\Gamma_n)})$ and $\big(\restr{r_n}{M_n}\big)^{-1}:f(\lambda^n B_{\ell_\infty(\Gamma_n)})\to M_n$ is $(3\lambda+2)$-Lipschitz.
\begin{proof}
Let us first prove that $r_n(M_n)=f(\lambda^n B_{\ell_\infty(\Gamma_n)})$. Clearly, the inclusion $r_n(M_n)\subset f(\lambda^n B_{\ell_\infty(\Gamma_n)})$ is satisfied. Let us prove the other inclusion. In fact, for $n=1$ it is trivially true and if $n\ge2$ then
$$\begin{aligned}r_n(M_n)=&r_n\Big( M_{n-1}\cup f\circ i_n\Big(f\big( \lambda^nB_{\ell_{\infty}(\Gamma_n)} \big)\setminus r_n\big( M_{n-1} \big)\Big) \Big)\\\supset& r_n\big(M_{n-1}\big)\cup\Big( r_n\circ f\circ i_n\circ f \big( \lambda^nB_{\ell_{\infty}(\Gamma_n)} \big)\setminus  r_n\circ f\circ i_n\circ r_n\big( M_{n-1} \big) \Big)\\=&r_n\big( M_{n-1} \big)\cup \Big(f\big(\lambda^n B_{\ell_\infty(\Gamma_n)}\big)\setminus r_n\big(M_{n-1}\big) \Big)=f\big(\lambda^n B_{\ell_\infty(\Gamma_n)}\big).\end{aligned}$$
We continue by proving that $\restr{r_n}{M_n}$ is injective. We proceed by induction in $n\in\N$. Clearly, $\restr{r_1}{M_1}$ is injective. Assuming that for arbitrary $n\ge2$ the retraction $\restr{r_{n-1}}{M_{n-1}}$ is injective, it suffices to prove that $\restr{r_n}{M_n}$ is also injective. Let us take $x,y\in M_n$ such that $r_n(x)=r_n(y)$. From the assumption it follows that $\restr{r_n}{M_{n-1}}$ is injective so we may restrict to the case when $x\notin M_{n-1}$. In this case, we claim that $y\notin M_{n-1}$. In fact, $r_n(y)=r_n(x)\in r_n(M_n\setminus M_{n-1})$ with
\begin{equation}\label{rnmn}\begin{aligned}r_n(M_n\setminus M_{n-1})=&r_n\big(f\circ i_n\big(f(\lambda^n B_{\ell_\infty(\Gamma_n)})\setminus r_n(M_{n-1}) \big)\big)\\=&r_n\big(f\circ i_n\big(r_n(M_n)\setminus r_n(M_{n-1}) \big)\big)\\=&r_n(M_n)\setminus r_n(M_{n-1}).\end{aligned}\end{equation}
This means that $r_n(y)\notin r_n(M_{n-1})$ so $y\notin M_{n-1}$. Then, it is enough to prove the case when $x,y\in M_n\setminus M_{n-1}$. In this case,
$$x=f\circ i_{n}\circ r_{n}(x)=f\circ i_{n}\circ r_{n}(y)=y,$$
and we have shown that $\restr{r_n}{M_n}$ is injective. Finally, let us prove that $\big(\restr{r_n}{M_n}\big)^{-1}$ is $(3\lambda+2)$-Lipschitz. We take two points $x,y\in M_n$ with $r_n(x)\neq r_n(y)$. Then, by Lemma \ref{dense} and the fact that $||r_n(x)-r_n(y)||\ge 1$,
$$\begin{aligned}\big|\big|\big(\restr{r_n}{B_n}\big)^{-1}( r_n(x))&-\big(\restr{r_n}{B_n}\big)^{-1}(r_n(y))\big|\big|=||x-y||\\\le&\big|\big|i_n\circ r_n(x)-x\big|\big|+\big|\big|i_n\circ r_n(y)-y\big|\big|+||i_n\circ r_n(x)-i_n\circ r_n(y)||\\\le& 2\lambda+2+\lambda||r_n(x)-r_n(y)||\le (3\lambda+2)||r_n(x)-r_n(y)||,\end{aligned}$$
which was to be proved.
\end{proof}
\end{lemma}

We may assume without loss of generality that $\lambda\in\N$. Take $s_n=\lambda^n$ for every $n\in\N$ and define the retraction $\phi_n:M\to M_{n}$ as $\phi_n=\big(\restr{r_n}{M_n}\big)^{-1}\circ T_{s_{n}}\circ r_{n}$. Let us establish some useful properties of these retractions.

\begin{prop}\label{globaldef}
For every $m,n\in\N$ it holds that $\phi_n\circ\phi_m=\phi_{\min(m,n)}$ and $||\phi_n||_{\text{Lip}}\le 3\lambda+2$.
\begin{proof}
The case when $n\ge m$ is clearly satisfied by the definition of retraction. Let us assume that $n<m$. It is easy to check that in this case $r_n\circ r_m=r_{n}$ and $\big(T_{s_n}\circ r_n\big)\circ\big(T_{s_m}\circ r_m\big)=T_{s_n}\circ r_n$. Hence,
$$\begin{aligned}\phi_n\circ \phi_m=&\big(\restr{r_n}{M_n}\big)^{-1}\circ T_{s_{n}}\circ r_{n}\circ \big(\restr{r_m}{M_m}\big)^{-1}\circ T_{s_{m}}\circ r_{m}\\=&\big(\restr{r_n}{M_n}\big)^{-1}\circ T_{s_{n}}\circ r_{n}\circ r_m\circ \big(\restr{r_m}{M_m}\big)^{-1}\circ T_{s_{m}}\circ r_{m}\\=&\big(\restr{r_n}{M_n}\big)^{-1}\circ\big( T_{s_{n}}\circ r_{n}\big)\circ\big( T_{s_{m}}\circ r_{m} \big)=\big(\restr{r_n}{M_n}\big)^{-1}\circ T_{s_{n}}\circ r_{n}=\phi_n. \end{aligned}$$
Finally, by Lemma \ref{Lipmain1} we conclude that $||\phi_n||_{\text{Lip}}\le 3\lambda+2$.
\end{proof}
\end{prop}

This finishes the first step. The only difference between the constructions made in this step and a retractional basis is that the number of points of $M_n\setminus M_{n-1}$ is not one in general. In the next steps we will sharpen this result by defining blocks between $M_n$ and $M_{n-1}$ whose difference consist of only one element, and the corresponding commutative Lipschitz retractions.

\subsection{Step 2}

In this step we are going to define intermediate blocks $D_n$ such that $M_n\subset D_n\subset M_{n+1}$ in order  to retract the elements of $M$ one by one.

To this end, we define an increasing sequence of sets $(C_n)$ as
$$C_n=f\circ i_{n+1}\circ f\circ T_{s_{n+1}}\circ r_{n+1}\circ i_n\big( f(s_{n+1}B_{\ell_\infty(\Gamma_n)})\setminus f( s_{n}B_{\ell_\infty(\Gamma_n)})\big),$$
for every $n\in\N$.
\begin{lemma}\label{contCB}
For every $n\in\N$, $C_n\subset M_{n+1}\setminus M_n$.
\begin{proof}
We denote by $A=f\circ T_{s_{n+1}}\circ r_{n+1}\circ i_n\big( f(s_{n+1}B_{\ell_\infty(\Gamma_n)})\setminus f( s_{n}B_{\ell_\infty(\Gamma_n)})\big)$
We claim that for every $n\in\N$,
\begin{equation}\label{cont1}\begin{aligned}A\subset r_{n+1}(M_{n+1})\setminus r_{n+1}(M_n).
\end{aligned}\end{equation}
Clearly $A\subset f(s_{n+1}B_{\ell_{\infty}(\Gamma_{n+1})})=r_{n+1}(M_{n+1})$ so it only remains  to prove that $A\cap r_{n+1}(M_n)=\emptyset$. Indeed, $r_n(A)=f(s_{n+1}B_{\ell_\infty(\Gamma_n)})\setminus f( s_{n}B_{\ell_\infty(\Gamma_n)})=f(s_{n+1}B_{\ell_\infty(\Gamma_n)})\setminus r_n(M_n)$ meaning that $r_n(A)\cap r_n(M_n)=\emptyset$. Hence,
$$r_n(A\cap r_{n+1}(M_n))\subset r_n(A)\cap r_n(M_n)=\emptyset,$$
so $A\cap r_{n+1}(M_n)=\emptyset$. This proves \eqref{cont1} and it follows that
$$C_n=f\circ i_{n+1}(A)\subset f\circ i_{n+1}\big(r_{n+1}(M_{n+1})\setminus r_{n+1}(M_n)\big)=M_{n+1}\setminus M_n.$$
\end{proof}
\end{lemma}

It now makes sense to define $D_n=M_n\cup C_n$ for every $n\in\N$. Then, we have the followig increasing chain of subsets of $M$
$$M_1\subset D_1\subset M_2\subset D_2\subset\cdots\subset M_n\subset D_n\subset\cdots$$
We are going to make use of the following Lemma \ref{Lipmain} which is an analogue of Lemma \ref{Lipmain1} in this new setting.

\begin{lemma}\label{Lipmain}
For every $n\in\N$ the map $\restr{r_n}{D_n}$ is injective. Moreover, $r_n(D_n)=r_n(M_{n+1})=f(s_{n+1}B_{\ell_\infty(\Gamma_n)})$ and $\big(\restr{r_n}{D_n}\big)^{-1}:f(s_{n+1}B_{\ell_\infty(\Gamma_n)})\to D_n$ is $(\lambda^2+2\lambda+2)$-Lipschitz.
\begin{proof}
Let us first show that $r_{n}(D_n)=r_n(M_{n+1})=f(s_{n+1}B_{\ell_\infty(\Gamma_n)})$. In fact,
$$\begin{aligned}r_n(D_n)=&r_n(M_n)\cup r_n(C_n)\\=&f(s_n B_{\ell_\infty(\Gamma_n)})\cup \big(f(s_{n+1} B_{\ell_\infty(\Gamma_n)})\setminus f (s_n B_{\ell_\infty(\Gamma_n)})\big)\\=&f(s_{n+1} B_{\ell_\infty(\Gamma_n)}).\end{aligned}$$
By Lemma \ref{Lipmain1},
$$\begin{aligned}r_n(M_{n+1})=&r_n\circ r_{n+1}(M_{n+1})=r_n\big( f(s_{n+1}B_{\ell_{\infty}(\Gamma_{n+1})}) \big)=f(s_{n+1}B_{\ell_\infty(\Gamma_n)}).\end{aligned}$$
We proceed by proving that $\restr{r_n}{D_n}$ is injective.  Let us take $x,y\in D_n$ such that $r_n(x)=r_n(y)$. We may assume that $x\notin M_n$ because the case when $x,y\in M_n$ was proven in Lemma \ref{Lipmain1}. As $x\in C_n$, by Lemma \ref{contCB} and the equality \eqref{rnmn} it follows that $r_n(y)=r_n(x)\in r_n(M_{n+1}\setminus M_n)=r_n(M_{n+1})\setminus r_n(M_n)$ so $y\in D_n\setminus M_n=C_n$. Therefore, $x,y\in C_n$ so that $x=f\circ i_{n+1}\circ f\circ T_{s_{n+1}}\circ  r_{n+1}\circ i_n(r_n(x))=f\circ i_{n+1}\circ f\circ T_{s_{n+1}}\circ  r_{n+1}\circ i_n(r_n(y))=y$, which proves that $\restr{r_n}{D_n}$ is injective.  Finally, to prove that $\Big|\Big|\Big(\restr{r_n}{D_{n}}\Big)^{-1}\Big|\Big|_{\text{Lip}}\le\lambda^2+2\lambda+2$ we take $x,y\in D_n$ with $r_n(x)\neq r_n(y)$. If $x,y\in M_{n}$ then by Lemma \ref{Lipmain1},
$$||x-y||\le (3\lambda+2)||r_n(x)-r_n(y)||.$$
Now, if $x\in D_{n}\setminus M_{n}$ and $y\in M_{n}$ then
$$x=f\circ i_{n+1}\circ f\circ T_{s_{n+1}}\circ  r_{n+1}\circ i_n(r_n(x)),$$
and by Lemma \ref{dense}, $||y-y_1||\le\lambda+1$ where $y_1=i_n\circ r_n(y)$. Clearly $||y_1||\le \lambda||r_n(y)||\le \lambda^{n+1}=s_{n+1}$ so that $T_{s_{n+1}}(y_1)=y_1$. Also, by the compatibility condition of the extensions we have that $y_1=i_{n+1}\circ r_{n+1}(y_1)$. Hence,
$$y_1=i_{n+1}\circ r_{n+1}\circ T_{s_{n+1}}\circ i_n\circ r_n(y)=i_{n+1}\circ T_{s_{n+1}}\circ r_{n+1}\circ i_n\circ r_n(y).$$
We let $y_2=i_{n+1}\circ f\circ T_{s_{n+1}}\circ r_{n+1}\circ i_n\circ r_n(y)$ so that $||y_2-y_1||\le \lambda$. This case is now done since
$$\begin{aligned}||x-y||\le&||x-f(y_2)||+||f(y_2)-y_2||+||y_2-y_1||+||y_1-y||\\\le& ||f\circ i_{n+1}\circ f\circ T_{s_{n+1}}\circ r_{n+1}\circ i_n||_{\text{Lip}}||r_n(x)-r_n(y)||+2\lambda+2\\\le& (\lambda^2+2\lambda+2)||r_n(x)-r_n(y)||.\end{aligned}$$
The last case is when $x,y\in D_n\setminus M_{n}$ so that
$$x=f\circ i_{n+1}\circ f\circ T_{s_{n+1}}\circ  r_{n+1}\circ i_n(r_n(x)),$$
$$y=f\circ i_{n+1}\circ f\circ T_{s_{n+1}}\circ  r_{n+1}\circ i_n(r_n(y)).$$
Finally, 
$$\begin{aligned}||x-y||\le& ||f\circ i_{n+1}\circ f\circ T_{s_{n+1}}\circ r_{n+1}\circ i_n||_{\text{Lip}}||r_n(x)-r_n(y)||\\\le&\lambda^2||r_n(x)-r_n(y)||.\end{aligned}$$
\end{proof}
\end{lemma}

Lemma \ref{Lipmain} allows us to define the intermediate retractions $\Psi_n:M_n\to D_{n-1}$ as $\Psi_n=\big(\restr{r_{n-1}}{D_{n-1}}\big)^{-1}\circ r_{n-1}$ for every $n\ge2$. Clearly, $\Psi_n\circ\Psi_m=\Psi_{\min(n,m)}$ for every $n,m\ge2$. These retractions factorize $\phi_{n-1}:M\to M_{n-1}$ making use of $D_{n-1}$ and $M_n$. More precisely,
$$\phi_{n-1}=\restr{\phi_{n-1}}{D_{n-1}}\circ \Psi_n\circ \phi_{n}.$$
In fact $\Psi_{n}(M_{n}\setminus M_{n-1})=C_n$ and $\restr{\Psi_{n}}{M_{n-1}}=Id_{M_{n-1}}$. The plan for the third step is to factorize $\Psi_n$  by certain retractions that only retract one element while in the forth step we will do the same for $\restr{\phi_{n-1}}{D_{n-1}}$.

\subsection{Step 3}

This is the longest and most technical step.

Our goal is to relabel the points of $M_n\setminus D_{n-1}$ as a finite sequence $(x^n_i)_{i=1}^{i(n)}$ and define Lipschitz retractions from $M_n$ onto $D_{n-1}\cup(x^n_l)_{l\le i}$ for every $i\in\{1,\dots, i(n)\}$ which factorize $\Psi_n$. More precisely, in this step we prove the following result.

\begin{prop}\label{step3}
For every $n\in\N$, $n\ge2$ there exists an order $(x_i^n)_{i=1}^{i(n)}=M_n\setminus D_{n-1}$ and for every $i\in\{0,\dots,i(n)-1\}$ there exists a retraction $\Psi_{n,i}:M_n\to D_{n-1}\cup(x^n_l)_{l\le i}$ such that for every $i,i_1,i_2\in\{0,\dots,i(n)-1\}$,
\begin{enumerate}
\item $\Psi_{n,i_1}\circ\Psi_{n,i_2}=\Psi_{n,\min(i_1,i_2)}$.
\item $\Psi_{n,0}=\Psi_n$.
\item $\Psi_{n,i}$ is $K$-Lipschitz for some $K$ independent of $i$.
\end{enumerate}
\end{prop}

Throughout this step we fix some $n\in\N$ with $n\ge2$.  We will make use of the following notation.

If $A,B\subset \Gamma$ are disjoint subsets, we define the operation $\oplus:\ell_\infty(A)\times \ell_\infty(B)\to\ell_\infty(A\cup B)$ as
$$x\oplus y(\gamma)=\begin{cases} x(\gamma)\;\;&\text{if }\gamma\in A,\\y(\gamma)&\text{if }\gamma\in B. \end{cases}$$
Given $C_A\subset \ell_\infty(A)$ and $C_B\subset \ell_\infty(B)$ we denote 
$$C_A\oplus C_B=\{x\oplus y\in\ell_\infty(A\cup B)\;:\;x\in C_A,\;y\in C_B\}.$$
With this terminology, we may identify $M_n$ with $r_n(M_n)$ by Lemma \ref{Lipmain1} and it is easy to see that $r_n(M_n)=r_{n-1}(D_{n-1})\oplus f(s_nB_{\ell_\infty(\Delta_n)})$ where $\Delta_n=\Gamma_n\setminus \Gamma_{n-1}$. We need to find a way to retract the points of $f(s_nB_{\ell_\infty(\Delta_n)})$ to an arbitrary $d\in \restr{D_{n-1}}{\Delta_n}$ one by one in a Lipschitz manner. Our approach is to first define a retractional basis in $f(\ell_\infty(\Delta_n))$ whose first element in the order is 0 and then translate its retractions by the element $d$.

Denote $N^n=f(\ell_\infty(\Delta_n))$ and for every $m\in\N\cup\{0\}$ we also denote $N^n_m=f\big(mB_{\ell_\infty(\Delta_n)}\big)$. It is clear that $N^n_m\setminus N^n_{m-1}=f\big(mS_{\ell_\infty(\Delta_n)}\big)$ for every $m\in\N$. Let us denote the points of $N^n_m\setminus N^n_{m-1}$ as $(y^{n}_{m,1},\dots y^{n}_{m,p_{n,m}})$ where $p_{n,m}=\#\big(N^n_m\setminus N^n_{m-1}\big)$ for every $m\in\N$. We consider  the set of indices  $A(n)=\bigcup\limits_{m\in\N}\bigcup\limits_{p=1}^{p_{n,m}}(m,p)$. Let us order the set $A(n)=(a^n_j)_{j\in\N}$ using the lexicographic order. Then, $N^n=(y^{n}_{a^n_j})_{j\in\N\cup\{0\}}$ where $y^n_{a^n_0}=0$. To simplify the notation we rename $y^n_j=y^n_{a^n_j}$ for every $j\in\N\cup\{0\}$. It is immediate that the order respects the norm, that is, if $j_1\le j_2\in \N\cup\{0\}$ then $||y^n_{j_1}||\le||y^n_{j_2}||$.

Let us denote $N^{n,j}=(y^n_{l})_{l\le j}$. Then, we may define the 'local' retractions of $N^n$, $T^n_j:N^{n,j}\to N^{n,j-1}$ as the retraction such that $T^n_j(y^n_{j})=T_{||y^n_{j}||-1}(y^n_{j})$ for every $j\in \N$. It is straightforward from the definition of retraction that
\begin{equation}\label{commuT}T^n_{j_1}\circ T^n_{j_2}=T^n_{j_2}\;\;\;\forall j_1\ge j_2.\end{equation}
Given a pair $j_1,j_2\in \N\cup\{0\}$ with $j_1<j_2$, we define the value $m(j_2,j_1)\in\N\cup\{0\}$ as the unique element of $\N\cup\{0\}$ such that
$$T_{j_1+1}^n\circ\cdots\circ T^n_{j_2}(y_{{j_2}}^n)=T_{m(j_2,j_1)}(y_{{j_2}}^n).$$
Let us point out that $m(j_2,j_1)$ is necessarily unique since $T_{j_1+1}^n\circ\cdots\circ T^n_{j_2}(y_{{j_2}}^n)\neq y_{j_2}^n$. We start the proof of Proposition \ref{step3} which will
be splitted into a series of auxiliary lemmas.

\begin{lemma}\label{bound}
Let $j_1,j_2\in\N\cup\{0\}$ be such that $j_1<j_2$. Then $||y^n_{j_1}||-1\le m(j_2,j_1)\le||y^n_{j_1}||$.
\begin{proof}
 Clearly, if $y\in N^{n,j}$ then
\begin{equation}\label{induc}||T^n_j(y)||\ge \min\{||y||,||y_j^n||-1\}.\end{equation}
Also, by induction in $j_2-j_1$ using \eqref{induc} it is easy to prove that for every $y\in N^{n,j_2}$,
$$||T_{j_1+1}^n\circ\cdots\circ T^n_{j_2}(y)||\ge\min\{||y||,||y^n_{j_1}||-1\}.$$
Then, it is straightforward to see that $m(j_2,j_1)=||T_{j_1+1}^n\circ\cdots\circ T^n_{j_2}(y^n_{j_2})||\ge \min\{||y^n_{j_2}||,||y^n_{j_1}||-1\}=||y^n_{j_1}||-1$. Finally, since $T_{j_1+1}^n\circ\cdots\circ T^n_{j_2}(N^{n,j_2})=N^{n,j_1}$ we know that $m(j_2,j_1)\le ||y^n_{j_1}||$.
\end{proof}
\end{lemma}

\begin{remark}
Clearly $m(j_2,j_1)$ depends on $n\in\N$ but since  $n\ge2$ is fixed throughout the step 3 and the value $m(j_2,j_1)$ is used only within this step we do not express the dependence on $n$.
\end{remark}

We are now ready to use the above constructions in our setting. Let us pick an arbitrary order in $D_{n-1}$ given by $(d^n_{1},\dots,d^n_{k(n)})$ where $k(n)=\#\big( D_{n-1} \big)$. For every $j\in\N\cup\{0\}$ we consider the sets of indices $K(j,n)=\{k\in\{1,\dots,k(n)\}\;:\; \restr{d^n_k}{\Delta_n}+y^n_{j}\in N^n_{s_n}\}$ and $J(n)=\{j\in\N\;:\;K(j,n)\neq\emptyset\}$. It is easy to see that $J(n)$ is finite. Finally, let us define the set of indices $E(n)=\{(j,k)\;:\; j\in J(n),\;k\in K(j,n)\}$, which is also a finite set. Now, we rename $E(n)$ as $\{e^n_1,\dots,e^n_{i(n)}\}$ following the lexicographic order. If $i\in \{1,\dots,i(n)\}$ with $e^n_i=(j,k)$ then we denote $e_{i}^n(1)=j$ and $e^n_i(2)=k$. The fact that $(e_i)_{i=1}^{i(n)}$ is ordered lexicographically may be understood as follows. $i_1,i_2\in\{1,\dots,i(n)\}$ with $i_1<i_2$ means that $e^n_{i_1}(1)\le e^n_{i_2}(1)$ and $e^n_{i_1}(1)=e^n_{i_2}(1)$ implies $e^n_{i_1}(2)<e^n_{i_2}(2)$.

\begin{lemma}\label{deforder}
The following equivalences hold:
\begin{itemize}
\item $x\in M_n\setminus D_{n-1}$ if, and only if, there is a unique $e^n_i=(j,k)\in E(n)$ such that $x=\big(\restr{r_n}{M_n}\big)^{-1}\big(r_{n-1}(d^n_k)\oplus \big(\restr{d_k^n}{\Delta_n}+y_{j}^n\big)\big)$.
\item $x\in D_{n-1}$ if and only if there is a unique $k\in \{1,\dots,k(n)\}=K(0,n)$ such that $x=\big(\restr{r_n}{M_n}\big)^{-1}\big(r_{n-1}(d^n_k)\oplus \big(\restr{d_k^n}{\Delta_n}+y_{0}^n\big)\big)=d_k^n$.
\end{itemize}
\begin{proof}
It follows from the fact that $\big(\restr{r_n}{M_n}\big)^{-1}:r_n(M_n)\to M_n$ is a bijection and the fact that $r_n(M_n)=r_{n-1}(D_{n-1})\oplus N_{s_n}^n$.
\end{proof}
\end{lemma}

Following Lemma \ref{deforder}, we are able to establish an order of $M_n\setminus D_{n-1}$ by $(x^n_{i})_{i=1}^{i(n)}$ where $x^n_{i}=\big(\restr{r_n}{M_n}\big)^{-1}\big(r_{n-1}(d^n_k)\oplus \big(\restr{d_k^n}{\Delta_n}+y_{j}^n\big)\big)$ whenever $e^n_i=(j,k)\in E(n)$. Now for every $i\in\{1,\dots,i(n)\}$ we consider the set $D_{n-1,i}=D_{n-1}\cup(x^n_{l})_{l\le i}$ and $D_{n-1,0}=D_{n-1}$. This allows us to define the local intermediate retractions $\psi_{n,i}:D_{n-1,i}\to D_{n-1,i-1}$ for every $i\in\{1,\dots i(n)\}$ by
$$\psi_{n,i}(x)=\begin{cases}\Big(\restr{r_n}{M_n}\Big)^{-1}\Big( r_{n-1}(d_k^n)\oplus \big(\restr{d_k^n}{\Delta_n}+T_{||y^n_j||-1}(y^n_{j})\big) \Big)\;\;&\text{if }x=x^n_i,\\x&\text{otherwise},\end{cases}$$
where $e^n_i=(j,k)$. Finally, for every $i\in \{0,\dots, i(n)-1\}$ we define the global intermediate retractions as 
$$\Psi_{n,i}=\psi_{n,{i+1}}\circ\cdots\circ \psi_{n,{i(n)}}:M_n\to D_{n-1,i}.$$
As we pointed out in \eqref{commuT}, in this setting we also want to mention that
\begin{equation}\label{commupsi}\psi_{n,i_1}\circ \psi_{n,i_2}=\psi_{n,i_2}\;\;\;\text{if }\; i_1\ge i_2\in\{1,\dots,i(n)\},\end{equation}
and
\begin{equation}\label{commuPsi}\Psi_{n,i_1}\circ \Psi_{n,i_2}=\Psi_{n,\min(i_1,i_2)}\;\;\;\text{if }\; i_1,i_2\in \{0,\dots,i(n)-1\}.\end{equation}

\begin{lemma}\label{lemmamain}
Let $i_1,i_2\in\{1,\dots,i(n)\}$ be such that $i_1<i_2$. Then
\begin{equation}\label{defret}\Psi_{n,i_1}(x^n_{{i_2}})=\Big(\restr{r_n}{M_n}\Big)^{-1}\Big( r_{n-1}(d_{k_2}^n)\oplus \Big(\restr{d_{k_2}^n}{\Delta_n}+T_{M(i_1,i_2)}(y^n_{j_2})\Big) \Big),\end{equation}
where $e^n_{i_1}=({j_1},k_1),e^n_{i_2}=({j_2},k_2)\in E(n)$ and
$$M(i_1,i_2)=\begin{cases}m({j_2},j_1)\;\;&\text{if }k_2\le k_1,\\m({j_2},j_1-1)&\text{if }k_2>k_1.\end{cases}$$
\begin{remark}\label{advise}
$M(i_1,i_2)$ is well defined whenever $i_1<i_2$ since in this case $j_1<j_2$ provided $k_1\ge k_2$. Also $j_1-1<j_2$ provided $k_1<k_2$.
\end{remark}
\begin{proof}
We may assume that $\psi_{n,i_2}(x^n_{i_2})\notin D_{n-1}$ because if $\psi_{n,i_2}(x^n_{i_2})\in D_{n-1}$ then $\Psi_{n,i_1}(x^n_{i_2})=\psi_{n,i_2}(x^n_{i_2})\in D_{n-1}$ and the statement is satisfied with $M(i_1,i_2)=0$. Throughout all the proof we will denote $i_3\in \{1,\dots,i_2-1\}$ such that $x^n_{i_3}=\psi_{n,i_2}(x^n_{i_2})$ and $e_{i_3}^n=({j_3},k_2)$ for some $j_3<j_2$.  It holds that
\begin{equation}\label{rest}T^n_{j_2}(y^n_{j_2})=y^n_{j_3}.\end{equation}

\textbf{Case 1, $k_2\le k_1$.} In this case we prove \eqref{defret} for $M(i_1,i_2)=m({{j_2}},j_1)$ by induction in $i_2-i_1\in\N$. Let us recall that, as mentioned in Remark \ref{advise}, $j_1<j_2$ throughout this case.

If $i_2-i_1=1$ then $i_3\le i_2-1=i_1$. Hence,
$$\begin{aligned}\Psi_{n,i_1}(x^n_{{i_2}})=&\psi_{n,i_2}(x^n_{i_2})=\Big(\restr{r_n}{B_n}\Big)^{-1}\Big( r_{n-1}(d_{k_2}^n)\oplus \Big(\restr{d_{k_2}^n}{\Delta_n}+T_{j_2}^n(y^n_{j_2})\Big) \Big).\end{aligned}$$
It is then enough to show that $T_{j_2}^n(y^n_{j_2})=T_{m(j_2,j_1)}(y^n_{j_2})$. Following the lexicographic order $j_3=e^n_{i_3}(1)\le e^n_{i_1}(1)=j_1<j_2$, equivalently $j_3<j_1+1\le j_2$. There are two possible cases, namely $j_2=j_1+1$ or $j_2-1\ge j_1+1$. If $j_2=j_1+1$ then by definition $T_{j_2}^n(y^n_{j_2})=T_{m(j_2,j_1)}(y^n_{j_2})$. Otherwise, if $j_2-1\ge j_1+1$ then
$$T^n_{j_2}(y^n_{j_2})=y^n_{j_3}=T_{j_1+1}^n\circ\cdots\circ T^n_{j_2-1}(y^n_{j_3})=T_{j_1+1}^n\circ\cdots\circ T_{j_2}^n(y^n_{j_2})=T_{m(j_2,j_1)}(y^n_{j_2}).$$
This proves the first step of the induction.

Now we go on and prove the inductive step. Let us assume that $i_2-i_1\ge2$ and that
\begin{equation}\label{defret2}\Psi_{n,i_1}(x^n_{{i}})=\Big(\restr{r_n}{B_n}\Big)^{-1}\Big( r_{n-1}(d_{k}^n)\oplus \Big(\restr{d_{k}^n}{\Delta_n}+T_{m(j,j_1)}(y^n_{j})\Big) \Big)\end{equation}
holds for every index $i\in\{i_1+1,\cdots,i_2-1\}$ with $e_i^n=(j,k)$ and $k\le k_1$. We now split the proof into two different cases, namely $i_1+1\le i_3$ and $i_1+1>i_3$. If $i_1+1\le i_3$, then we claim that $j_1<j_3<j_2$. In fact, $i_1<i_3$ and $e^n_{i_1}(2)=k_1\ge k_2=e^n_{i_3}(2)$ implies that $j_1=e^n_{i_1}(1)<e^n_{i_3}(1)=j_3$. Also, $i_3<i_2$ and $e^n_{i_3}(2)=k_2=e^n_{i_2}(2)$ implies that $j_3=e^n_{i_3}(1)<e^n_{i_2}(1)=j_2$ so we have shown that $j_1<j_3<j_2$. Therefore, by \eqref{commuT},
$$\begin{aligned}T_{m( j_3,j_1)}(y^n_{j_3})=&T_{j_1+1}^n\circ\cdots\circ T^n_{j_3}\big(y^n_{j_3}\big)\\=&T_{j_1+1}^n\circ\cdots\circ T^n_{j_3}\big(T^n_{j_3+1}\circ\cdots\circ T^n_{j_2+1}(y^n_{j_3})\big)\\=&T_{j_1+1}^n\circ\cdots\circ T^n_{j_2}\big(y^n_{j_2}\big)=T_{m(j_2,j_1)}(y^n_{j_2}).\end{aligned}$$
Hence, by equality \eqref{commupsi} and the assumption \eqref{defret2} with $i=i_3$ we have that
$$\begin{aligned}\Psi_{n,i_1}(x^n_{{i_2}})=&\Psi_{n,i_1}(x^n_{i_3})=\Big(\restr{r_n}{B_n}\Big)^{-1}\Big( r_{n-1}(d_{k_2}^n)\oplus \Big(\restr{d_{k_2}^n}{\Delta_n}+T_{m( j_3,j_1)}(y^n_{j_3})\Big) \Big)\\=&\Big(\restr{r_n}{B_n}\Big)^{-1}\Big( r_{n-1}(d_{k_2}^n)\oplus \Big(\restr{d_{k_2}^n}{\Delta_n}+T_{m(j_2,j_1)}(y^n_{j_2})\Big) \Big).\end{aligned}$$
Otherwise, if $i_1+1>i_3$ then $i_2-1\ge i_1+1>i_3$ so by \eqref{commupsi},
$$\begin{aligned}\Psi_{n,i_1}(x^n_{i_2})=&\psi_{n,i_1+1}\circ \cdots\circ \psi_{n,i_2-1}(x^n_{i_3})=x^n_{i_3}\\=&\Big(\restr{r_n}{B_n}\Big)^{-1}\Big( r_{n-1}(d_{k_2}^n)\oplus \Big(\restr{d_{k_2}^n}{\Delta_n}+T_{j_2}^n(y^n_{j_2})\Big) \Big).\end{aligned}$$
It only remains to see that $T_{j_2}^n(y^n_{j_2})=T_{m({j_2},j_1)}(y^n_{j_2})$. Since $i_3\le i_1$ then $j_3\le j_1<j_2$. If $j_1+1=j_2$ then $T_{j_2}^n(y^n_{j_2})=T_{m({j_2},j_1)}(y^n_{j_2})$ by definition. Otherwise, $j_3<j_1+1\le j_2-1$ so we finish the proof with
$$T_{j_2}^n(y^n_{j_2})=y^n_{j_3}=T_{j_1+1}^n\circ \cdots\circ T_{j_2-1}^n(y^n_{j_3})=T_{j_1+1}^n\circ \cdots\circ T_{j_2}^n(y^n_{j_2})=T_{m({j_2},j_1)}(y^n_{j_2}).$$

\textbf{Case 2, $k_2>k_1$.} In this case we prove \eqref{defret} with $M(i_1,i_2)=m(j_2,j_1-1)$. Again, let us prove this case by induction in $i_2-i_1\in\N$. If $i_2-i_1=1$ then
$$\begin{aligned}\Psi_{n,i_1}(x^n_{{i_2}})=&\psi_{n,i_2}(x^n_{i_2})=\Big(\restr{r_n}{B_n}\Big)^{-1}\Big( r_{n-1}(d_{k_2}^n)\oplus \Big(\restr{d_{k_2}^n}{\Delta_n}+T_{j_2}^n(y^n_{j_2})\Big) \Big).\end{aligned}$$
Thus, it is enough to show that $T^n_{j_2}(y^n_{j_2})=T_{m(j_2,j_1-1)}(y^n_{j_2})$. Since $i_3\le i_2-1=i_1$ with $e^n_{i_3}(2)=k_2>k_1=e^n_{i_1}(2)$ we deduce that  $j_3=e^n_{i_3}(1)< e^n_{i_1}(1)=j_1$. That is, $j_3<j_1\le j_2$ and we distinguish two cases, namely $j_1=j_2$ and $j_1\le j_2-1$. If $j_1=j_2$ then $T^n_{j_2}(y^n_{j_2})=T_{m(j_2,j_1-1)}(y^n_{j_2})$ by definition. Otherwise $j_3<j_1\le j_2-1$ and therefore,
$$T^n_{j_2}(y^n_{j_2})=y^n_{j_3}=T_{j_1}^n\circ\cdots\circ T^n_{j_2-1}(y^n_{j_3})=T_{j_1}^n\circ\cdots\circ T_{j_2}^n(y^n_{j_2})=T_{m(j_2,j_1-1)}(y^n_{j_2}).$$
Now that we have proven the first step of the induction we pass to the inductive step when $i_2-i_1\ge2$. We assume here that
\begin{equation}\label{defret3}\Psi_{n,i_1}(x^n_{{i}})=\Big(\restr{r_n}{B_n}\Big)^{-1}\Big( r_{n-1}(d_{k}^n)\oplus \Big(\restr{d_{k}^n}{\Delta_n}+T_{m(j,j_1-1)}(y^n_{j})\Big) \Big)\end{equation}
holds for every index $i\in\{i_1+1,\dots,i_2-1\}$ with $e_i^n=(j,k)\in E(n)$. We again split the proof of the inductive step into two different cases, namely, $i_1+1\le i_3$ and $i_1+1>i_3$. If $i_1+1\le i_3$ then $j_1\le j_3< j_2$ so by \eqref{commuT},
$$\begin{aligned}T_{m( j_3,j_1-1)}(y^n_{j_3})=&T_{j_1}^n\circ\cdots\circ T^n_{j_3}(y^n_{j_3})=T_{j_1}^n\circ\cdots\circ T^n_{j_3}\big(T^n_{j_3+1}\circ\cdots\circ T^n_{j_2}(y^n_{j_3})\big)\\=&T_{j_1}^n\circ\cdots\circ T^n_{j_2}(y^n_{j_2})=T_{m(j_2,j_1-1)}(y^n_{j_2}).\end{aligned}$$
Hence, using \eqref{defret3} with $i=i_3$,
$$\begin{aligned}\Psi_{n,i_1}(x^n_{{i_2}})=&\Psi_{n,i_1}(x^n_{i_3})=\Big(\restr{r_n}{B_n}\Big)^{-1}\Big( r_{n-1}(d_{k_2}^n)\oplus \Big(\restr{d_{k_2}^n}{\Delta_n}+T_{m(j_3,j_1-1)}(y^n_{j_3})\Big) \Big)\\=&\Big(\restr{r_n}{B_n}\Big)^{-1}\Big( r_{n-1}(d_{k_2}^n)\oplus \Big(\restr{d_{k_2}^n}{\Delta_n}+T_{m(j_2,j_1-1)}(y^n_{j_2})\Big) \Big).\end{aligned}$$
It only remains to treat the case when $i_1+1>i_3$, equivalently $i_3\le i_1$. We first claim that
\begin{equation}\label{eqinduc1}T_{j_2}^n(y^n_{j_2})=T_{m(j_2,j_1-1)}(y^n_{j_2}).\end{equation}
Since $e^n_{i_3}(2)=k_2>k_1=e_{i_1}^n(2)$, by means of the lexicographic order it must hold that $j_3=e^n_{i_3}(1)<e^n_{i_1}(1)=j_1\le j_2$. We split  the proof of \eqref{eqinduc1} into two different cases, $j_1=j_2$ and $j_1\le j_2-1$. If $j_1=j_2$ then \eqref{eqinduc1} holds by definition. Otherwise, $j_3<j_1\le j_2-1$ so that
$$\begin{aligned}T^n_{j_2}(y^n_{j_2})=y^n_{j_3}=&T^n_{j_1}\circ \cdots\circ T^n_{j_2-1}(y^n_{j_3})=T_{j_1}^n\circ\cdots\circ T_{j_2}^n(y^n_{j_2})=T_{m({j_2},j_1-1)}(y^n_{j_2}).\end{aligned}$$
This finishes the proof since using \eqref{commupsi} it follows that
$$\begin{aligned}\Psi_{n,i_1}(x^n_{i_2})=&\Psi_{n,i_1}(x^n_{i_3})=\psi_{n,i_1+1}\circ \cdots\circ \psi_{n,i_2-1}(x^n_{i_3})=x^n_{i_3}\\=&\Big(\restr{r_n}{B_n}\Big)^{-1}\Big( r_{n-1}(d_{k_2}^n)\oplus \Big(\restr{d_{k_2}^n}{\Delta_n}+T_{j_2}^n(y^n_{j_2})\Big) \Big)\\=&\Big(\restr{r_n}{B_n}\Big)^{-1}\Big( r_{n-1}(d_{k_2}^n)\oplus \Big(\restr{d_{k_2}^n}{\Delta_n}+T_{m({j_2},j_1-1)}(y^n_{j_2})\Big) \Big).\end{aligned}$$
\end{proof}
\end{lemma}

\begin{lemma}\label{remlip}
For every $i_1<i_2\in\{1,\dots,i(n)\}$,
$$||y^n_{j_1}||-2\le M(i_1,i_2)\le ||y^n_{j_1}||,$$
where $j_1=e^n_{i_1}(1)$ and $j_2=e^n_{i_2}(1)$.
\begin{proof}
By Lemma \ref{bound} it follows that
$$M(i_1,i_2)\le\max\{m(j_2,j_1),m(j_2,j_1-1)\}\le\max\{||y^n_{j_1}||,||y^n_{j_1-1}||\}=||y^n_{j_1}||.$$
Also, it is straightforward to check that $||y^n_{j_1-1}||\ge||y^n_{j_1}||-1$ so again using Lemma \ref{bound} we obtain that
$$M(i_1,i_2)\ge\min\{m(j_2,j_1),m(j_2,j_1-1)\}\ge\min\{||y^n_{j_1}||-1,||y^n_{j_1-1}||-1\}\ge||y^n_{j_1}||-2.$$
\end{proof}
\end{lemma}

\begin{prop}
It holds that $\Psi_{n,0}=\Psi_n$.
\begin{proof}
Let us take $i\in\{1,\dots,i(n)\}$. We are denoting for this proof $e_i^n=(j,k)$ and $e_1^n=(j_1,k_1)$. We claim that $||y^n_{j_1}||=1$.

As $x_1^n\notin D_{n-1}$ we know that $||y^n_{j_1}||\neq0$ so we only have to show that $||y^n_{j_1}||\le1$. If $\restr{d_{k_1}^n}{\Delta_n}=0$ then clearly $\restr{d_{k_1}^n}{\Delta_n}+y^n_1=y^n_1\in f(B_{\ell_\infty(\Delta_n)})\subset N^n_{s_n}$ so $1=\min\{j\in J(n)\;:\;k_1\in K(j,n)\}=j_1$. Hence, $||y^n_{j_1}||=||y^n_1||=1$. Otherwise, if $||\restr{d_{k_1}^n}{\Delta_n}||\ge1$ then clearly $\restr{d_{k_1}^n}{\Delta_n}-\frac{\restr{d_{k_1}^n}{\Delta_n}}{||\restr{d_{k_1}^n}{\Delta_n}||}\in N^n_{s_n}$. This means that there must exist a $j_2\in J(n)$ such that $y^n_{j_2}=-\frac{\restr{d_{k_1}^n}{\Delta_n}}{||\restr{d_{k_1}^n}{\Delta_n}||}$ and $k_1\in K(j_2,n)$. Hence, $j_2\ge \min\{j\in J(n)\;:\;k_1\in K(j,n)\}=j_1$ so $||y^n_{j_1}||\le||y^n_{j_2}||=1$ and the claim is proven.

Now, we pass to the proof that $\Psi_{n,0}=\Psi_n$ distinguishing two different cases. The first case is when $\Psi_{n,1}(x_i^n)=x_1^n$. In this case,
$$\begin{aligned}\Psi_{n,0}(x_i^n)=&\psi_{n,1}\circ\Psi_{n,1}(x_i^n)=\psi_{n,1}(x_1^n)\\=&\Big(\restr{r_n}{M_n}\Big)^{-1}\Big( r_{n-1}(d_{k_1}^n)\oplus \Big(\restr{d_{k_1}^n}{\Delta_n}+T_{||y^n_{j_1}||-1}(y^n_{j_1})\Big) \Big)\\=&\Big(\restr{r_n}{M_n}\Big)^{-1}\Big( r_{n}(d_{k}^n)\Big)\\=&d_k^n=\Big(\restr{r_{n-1}}{D_{n-1}}\Big)^{-1}\Big( r_{n-1}(x_i^n)\Big)=\Psi_n(x_i^n).\end{aligned}$$
The second case is when $\Psi_{n,1}(x_i^n)\neq x_1^n$ or equivalently $\Psi_{n,1}(x_i^n)\in D_{n-1,1}\setminus\{x_1^n\}=D_{n-1}$. In this case $i>1$ so by Lemma \ref{lemmamain},
$$\begin{aligned}r_{n-1}\big(\Psi_{n,1}(x_i^n)\big)=&r_{n-1}\bigg(\Big(\restr{r_n}{M_n}\Big)^{-1}\Big( r_{n-1}(d_{k}^n)\oplus \Big(\restr{d_{k}^n}{\Delta_n}+T_{M(1,i)}(y^n_{j})\Big) \Big)\bigg)\\=&r_{n-1}(d_k^n).\end{aligned}$$
Finally, since $\restr{r_{n-1}}{D_{n-1}}$ is an injection and $\Psi_{n,1}(x_i^n)\in D_{n-1}$ we have that
$$\begin{aligned}\Psi_{n,0}(x_i^n)=&\psi_{n,1}\circ\Psi_{n,1}(x_i^n)=\Psi_{n,1}(x_i^n)\\=&\Big(\restr{r_{n-1}}{D_{n-1}}\Big)^{-1}\Big(r_{n-1}\big(\Psi_{n,1}(x_i^n)\big)\Big)\\=&\Big(\restr{r_{n-1}}{D_{n-1}}\Big)^{-1}\Big(r_{n-1}\big(d_k^n\big)\Big)\\=&d_k^n=\Big(\restr{r_{n-1}}{D_{n-1}}\Big)^{-1}\Big( r_{n-1}(x_i^n)\Big)=\Psi_n(x_i^n).\end{aligned}$$
\end{proof}
\end{prop}

\begin{prop}
For every $i\in\{1,\dots,i(n)\}$ the retraction $\Psi_{n,i}$ is $(3\lambda+2)(3\lambda^2+6\lambda+11)$-Lipschitz.
\begin{proof}
Let us take distinct points $x_1,x_2\in M_n$, i.e. $||x_1-x_2||\ge1$. We may assume without loss of generality that  $x_1\in M_n\setminus D_{n-1,i}$ because otherwise $\Psi_{n,i}$ acts as an identity. Hence, there is $i_1\in\{i+1,\dots,i(n)\}$ with $e^n_{i_1}=({j_1},k_1)\in E(n)$ such that $x_1=x_{i_1}^n$. There are two possibilities, namely either $x_2\in D_{n-1,i}$ or $x_2\in M_n\setminus D_{n-1,i}$. If $x_2\in D_{n-1,i}$ then by Lemma \ref{deforder} there are unique $j_2\in J(n)\cup\{0\}$ and $k_2\in K(j_2,n)$ such that $({j_2},k_2)\le(j,k)=e^n_i$ (considering the lexicographic order) and $x_2=\big(\restr{r_n}{M_n}\big)^{-1}\big(r_{n-1}(d^n_{k_2})\oplus \big(\restr{d_{k_2}^n}{\Delta_n}+y_{{j_2}}^n\big)\big)$. As $j_2\le j\le j_1$, by Lemma \ref{remlip} we have $||y^n_{j_2}||\le||y^n_{j_1}||\le M(i,i_1)+2$ and so $\max\{0,||y^n_{j_2}||-M(i,i_1)\}\le 2$. Hence,
$$\begin{aligned}|| T_{M(i,i_1)}(y^n_{j_1})-y^n_{j_2} ||\le&|| T_{M(i,i_1)}(y^n_{j_1})-T_{M(i,i_1)}(y^n_{j_2})||+||T_{M(i,i_1)}(y^n_{j_2})-y^n_{j_2}||\\=&||T_{M(i,i_1)}(y^n_{j_1})-T_{M(i,i_1)}(y^n_{j_2})||+\max\{0,||y^n_{j_2}||-M(i,i_1)\}\\\le&||y^n_{j_1}-y^n_{j_1}||+2.\end{aligned}$$
Also, by Lemma \ref{bound}, it follows that
$$\begin{aligned}||\restr{d_{k_1}^n}{\Delta_n}-\restr{d_{k_2}^n}{\Delta_n}||\le& ||d^n_{k_1}-d^n_{k_2}||\\=&\Big|\Big|\Big(\restr{r_{n-1}}{D_{n-1}}\Big)^{-1}(r_{n-1}(d^n_{k_1}))-\Big(\restr{r_{n-1}}{D_{n-1}}\Big)^{-1}(r_{n-1}(d^n_{k_2}))\Big|\Big|\\=&\Big|\Big|\Big(\restr{r_{n-1}}{D_{n-1}}\Big)^{-1}(r_{n-1}(x_1))-\Big(\restr{r_{n-1}}{D_{n-1}}\Big)^{-1}(r_{n-1}(x_2))\Big|\Big|\\\le&(\lambda^2+2\lambda+2)||x_1-x_2||.\end{aligned}$$
Therefore,
$$\begin{aligned}||y^n_{j_1}-y^n_{j_2}||\le&\big|\big|\big( \restr{d^n_{k_1}}{\Delta_n}+y^n_{j_1} \big)-\big( \restr{d^n_{k_2}}{\Delta_n}+y^n_{j_2} \big)\big|\big|+\big|\big| \restr{d^n_{k_1}}{\Delta_n}-\restr{d^n_{k_2}}{\Delta_n} \big|\big|\\=& ||\restr{(x_1-x_2)}{\Delta_n}||+\big|\big| \restr{d^n_{k_1}}{\Delta_n}-\restr{d^n_{k_2}}{\Delta_n} \big|\big|\le (\lambda^2+2\lambda+3)||x_1-x_2||.\end{aligned}$$
By Lemma \ref{Lipmain1} and  Lemma \ref{lemmamain} we compute
$$\begin{aligned}||\Psi_{n,i}(x_1)-\Psi_{n,i}(x_2)||\le& \Big|\Big| \Big( \restr{r_n}{M_n} \Big)^{-1} \Big|\Big|_{\text{Lip}}\Big(|| r_{n-1}(d^n_{k_1}-d^n_{k_2}) ||\\&+||\restr{d_{k_1}^n}{\Delta_n}-\restr{d_{k_2}^n}{\Delta_n}||+|| T_{M(i,i_1)}(y^n_{j_1})-y^n_{j_2} ||\Big)\\&\le (3\lambda+2)(2\lambda^2+4\lambda+8)||x_1-x_2||.\end{aligned}$$
Finally, if $x_2\in M_n\setminus D_{n-1,i}$ then $x_2=x^n_{i_2}$, for some $i_2\in\{i+1,\dots,i(n)\}$. If we denote $e_{i_2}^n=(j_2,k_2)$, then by Lemma \ref{remlip}, 
$$|M(i,i_1)-M(i,i_2)|\le\big|\,||y^n_{j_1}||-||y^n_{j_2}||\,\big|+2\le||y^n_{j_1}-y^n_{j_2}||+2.$$
Therefore,
$$\begin{aligned}|| T_{M(i,i_1)}(y^n_{j_1})-T_{M(i,i_2)}(y^n_{j_2}) ||\le&||T_{M(i,i_1)}(y^n_{j_1})-T_{M(i,i_1)}(y^n_{j_2})||\\&+||T_{M(i,i_1)}(y^n_{j_2})-T_{M(i,i_2)}(y^n_{j_2})||\\\le&|y^n_{j_1}-y^n_{j_2}||+|M(i,i_1)-M(i,i_2)|\\\le&2||y^n_{j_1}-y^n_{j_2}||+2.\end{aligned}$$
Following the same computation as in the previous case we have that
$$||\restr{d_{k_1}^n}{\Delta_n}-\restr{d_{k_2}^n}{\Delta_n}||\le(\lambda^2+2\lambda+2)||x_1-x_2||,$$
and
$$||y^n_{j_1}-y^n_{j_2}||\le (\lambda^2+2\lambda+3)||x_1-x_2||.$$ 
By Lemma \ref{Lipmain1} and Lemma \ref{lemmamain} we finally conclude that
$$\begin{aligned}||\Psi_{n,i}(x_1)-\Psi_{n,i}(x_2)||\le& \Big|\Big| \Big( \restr{r_n}{M_n} \Big)^{-1} \Big|\Big|_{\text{Lip}}\Big(|| r_{n-1}(d^n_{k_1}-d^n_{k_2}) ||\\&+||\restr{d_{k_1}^n}{\Delta_n}-\restr{d_{k_2}^n}{\Delta_n}||+|| T_{M(i,i_1)}(y^n_{j_1})-T_{M(i,i_2)}(y^n_{j_2}) ||\Big)\\&\le (3\lambda+2)(3\lambda^2+6\lambda+11)||x_1-x_2||.\end{aligned}$$
\end{proof}
\end{prop}

This finishes the proof of Proposition \ref{step3} with $K=(3\lambda+2)(3\lambda^2+6\lambda+11)$.

\subsection{Step 4}
 
In this step we focus on the points that are in $D_n\setminus M_{n}=C_n$. Similarly to what we have done in the third step, we will prove the following Proposition \ref{step4}.

\begin{prop}\label{step4}
For every $n\in\N$, there exists an order $(c_i^n)_{i=1}^{c(n)}=D_n\setminus M_n=C_n$ and for every $i\in\{0,\dots,c(n)-1\}$ there exists a retraction $\phi_{n,i}:D_n\to M_{n}\cup(c^n_l)_{l\le i}$ such that for every $i,i_1,i_2\in\{0,\dots,c(n)-1\}$, we have
\begin{enumerate}
\item \label{1step4}$\phi_{n,i_1}\circ\phi_{n,i_2}=\phi_{n,\min(i_1,i_2)}$.
\item \label{2step4}$\phi_{n,0}=\restr{\phi_n}{D_n}$.
\item \label{3step4}$\phi_{n,i}$ is $L$-Lipschitz for some $L$ independent of $i$.
\end{enumerate}
\end{prop}

Again, let us fix an $n\in\N$. Thanks to Lemma \ref{Lipmain} we may identify $C_n$ with $r_n(C_n)=f(s_{n+1}B_{\ell_\infty(\Gamma_n)})\setminus f(s_{n}B_{\ell_\infty(\Gamma_n)})=\bigcup\limits_{k=s_n+1}^{s_{n+1}}f(kS_{\ell_\infty(\Gamma_n)})$. Now, it is a matter of ordering $\bigcup\limits_{k=s_n+1}^{s_{n+1}}f(kS_{\ell_\infty(\Gamma_n)})$. To do so, we give an arbitrary order to $f(kS_{\ell_\infty(\Gamma_n)})=N_k$ for $k\in\N\cap[s_n+1,s_{n+1}]$ as $N_k=(z_{k,1},\dots,z_{k,l_{n,k}})$. Here we consider the set of indices $G(n)=\bigcup\limits_{k=s_n+1}^{s_{n+1}}\bigcup\limits_{l=1}^{l_{n,k}}(k,l)$. We rename the points of $G(n)$ using the lexicographic order by $G(n)=(g^n_i)_{i=1}^{c(n)}$ where $c(n)=\#G(n)=\# C_n$. Let us denote $z_i^n=z_{g^n_i}$ for every $i\in\{1,\dots,c(n)\}$. The key 
property of this order is that  if $i_1,i_2\in\{1,\dots,c(n)\}$ with $i_1\le i_2$, then $||z^n_{i_1}||\le||z^n_{i_2}||$. This allows us to list the elements of $C_n$ as $c_{i}^n=\big(\restr{r_n}{C_n}\big)^{-1}(z^n_i)$. Considering $C_n=(c_{i}^n)_{i=1}^{c(n)}$ as an ordered set we obtain
$$||r_n(c_{i_1}^n)||=||z_{i_1}^n||\le||z_{i_2}^n||=||r_n(c_{i_2}^n)||\;\;\text{ whenever }i_1\le i_2\in\{1,\dots,c(n)\}.$$
Let $M_{n,i}=M_n\cup(c_l^n)_{l\le i}$ for $i\in\{1,\dots,c(n)\}$ and $M_{n,0}=M_n$. Then for every $i\in\{1,\dots,c(n)\}$ we define the local retraction $T_{n,i}:M_{n,i}\to M_{n,i-1}$ by
$$T_{n,i}(x)=\begin{cases}\big(\restr{r_n}{D_n}\big)^{-1}\big(T_{||r_n(x)||-1}(r_n(x))\big)\;\;&\text{if }x=c_i^n,\\x&\text{otherwise.}\end{cases}$$
We are finally ready to define the retractions of Proposition \ref{step4}. We define for each $i\in \{0,\dots,c(n)-1\}$ the retraction $\phi_{n,i}:D_n\to M_{n,i}$ by 
$$\phi_{n,i}=T_{n,i+1}\circ\cdots\circ T_{n,c(n)}.$$

\begin{lemma}\label{defstep4}
For every $i\in\{0,\dots,c(n)-1\}$ and $x\in D_n$,
$$\phi_{n,i}(x)=\big(\restr{r_n}{D_n}\big)^{-1}\big(T_{||r_n\circ\phi_{n,i}(x)||}(r_n(x))\big).$$
\begin{proof}
We proceed by induction in $c(n)-i\in\N$. In case when $c(n)-i=1$ we distinguish two cases. If $x=c_{c(n)}^n$ then
$$\phi_{n,c(n)-1}(x)=T_{n,c(n)}(x)=\big(\restr{r_n}{D_n}\big)^{-1}\big(T_{||r_n(x)||-1}(r_n(x))\big),$$
so that we only need to show that $||r_n(x)||-1=||r_n\circ\phi_{n,c(n)-1}(x)||$. In fact,
$$||r_n\circ\phi_{n,c(n)-1}(x)||=||T_{||r_n(x)||-1}(r_n(x))||=||r_n(x)||-1.$$
If $x\neq c_{c(n)}^n$ then
$$\phi_{n,c(n)-1}(x)=T_{n,c(n)}(x)=x=\big(\restr{r_n}{D_n}\big)^{-1}\big(T_{||r_n(x)||}(r_n(x))\big),$$
which ends the first step of the induction since $||r_n(x)||=||r_n\circ \phi_{n,c(n)-1}(x)||$. Finally, let us prove the induction step. If $c(n)-i\ge2$ we assume that
$$\phi_{n,i+1}(x)=\big(\restr{r_n}{D_n}\big)^{-1}\big(T_{||r_n\circ\phi_{n,i+1}(x)||}(r_n(x))\big).$$
Again, we just split this step into two cases. If $\phi_{n,i+1}(x)=c_{i+1}^n$ then
$$\begin{aligned}\phi_{n,i}(x)=&T_{n,i+1}\circ \phi_{n,i+1}(x)=\big(\restr{r_n}{D_n}\big)^{-1}\big(T_{||r_n\circ \phi_{n,i+1}(x)||-1}(r_n\circ \phi_{n,i+1}(x))\big)\\=&\big(\restr{r_n}{D_n}\big)^{-1}\big(T_{||r_n\circ \phi_{n,i+1}(x)||-1}(T_{||r_n\circ\phi_{n,i+1}(x)||}(r_n(x)))\big)\\=&\big(\restr{r_n}{D_n}\big)^{-1}\big(T_{||r_n\circ \phi_{n,i+1}(x)||-1}(r_n(x))\big)\\=&\big(\restr{r_n}{D_n}\big)^{-1}\big(T_{||r_n\circ \phi_{n,i}(x)||}(r_n(x))\big).\end{aligned}$$
Otherwise, if $\phi_{n,i+1}(x)\neq c_{i+1}^n$ then
$$\begin{aligned}\phi_{n,i}(x)=&T_{n,i+1}\circ \phi_{n,i+1}(x)=\phi_{n,i+1}(x)=\big(\restr{r_n}{D_n}\big)^{-1}\big(T_{||r_n\circ\phi_{n,i+1}(x)||}(r_n(x))\big)\\=&\big(\restr{r_n}{D_n}\big)^{-1}\big(T_{||r_n\circ\phi_{n,i}(x)||}(r_n(x))\big).\end{aligned}$$
\end{proof}
\end{lemma}

We are finally ready to prove the properties stated in Proposition \ref{step4}. Property $(\ref{1step4})$ of Proposition \ref{step4} follows from the definition of $\phi_{n,i}$. We prove properties $(\ref{2step4})$ and $(\ref{3step4})$ separately in the next Propositions.

\begin{prop}
The equality $\phi_{n,0}=\restr{\phi_n}{D_n}$ holds true.
\begin{proof}
If $x\in M_n$ then $\phi_{n,0}(x)=x=\phi_n(x)$. Let us take $c_i^n\in C_n$. Clearly, if some $x\in M_{n,i}$ is such that $||r_n(x)||\ge s_n$ then $||r_n\circ T_{n,i}(x)||\ge s_n$ so that iductively it follows that $||r_n\circ\phi_{n,0}(c_i^n)||=||r_n\circ T_{n,1}\circ\cdots\circ T_{n,c(n)}(c_i^n)||\ge s_n$. Also, we know that $r_n\circ \phi_{n,0}(c_i^n)\in r_n(M_n)=f(s_nB_{\ell_\infty(\Gamma_n)})$ so necessarily $||r_n\circ \phi_{n,0}(c_i^n)||=s_n$. We are done since Lemma \ref{defstep4} yields that
$$\begin{aligned}\phi_{n,0}(c_i^n)=&\big(\restr{r_n}{D_n}\big)^{-1}\big(T_{||r_n\circ\phi_{n,0}(c_i^n) ||}(r_n(c_i^n))\big)\\=&\big(\restr{r_n}{D_n}\big)^{-1}\big(T_{s_n}(r_n(c_i^n))\big)=\big(\restr{r_n}{M_n}\big)^{-1}\big(T_{s_n}(r_n(c_i^n))\big)=\phi_{n}(c_i^n).\end{aligned}$$
\end{proof}
\end{prop}

\begin{prop}
For every $i\in\{1,\dots,c(n)-1\}$, the retraction $\phi_{n,i}$ is $(2\lambda^2+4\lambda+4)$-Lipschitz.
\begin{proof}
Clearly, if $i\in \{0,\dots,c(n)-1\}$ and $x\in M_{n,i+1}$ then
\begin{equation}\label{ineqind}||r_n\circ T_{n,i+1}(x)||\ge \min\{ ||r_n(x)||,||z_{i}^n||-1 \}.\end{equation}
We claim that for every $x\in D_{n}$ and every $i\in\{0,\dots,c(n)-1\}$,
\begin{equation}\label{step2}||z^n_i||\ge||r_n\circ \phi_{n,i}(x)||\ge \min\{ ||r_n(x)||,||z^n_i||-1 \}.\end{equation}
The first inequality follows from the fact that $r_n\circ \phi_{n,i}(x)=z^n_{\widetilde i}\in r_n(M_{n,i})$ for some $\widetilde i\le i$. To prove the second inequality, we distinguish two cases.
If $x\in M_{n,i}$ then \eqref{step2} is satisfied since $r_n\circ \phi_{n,i}(x)=r_n(x)$ and clearly $||r_n(x)||\ge\min\{||r_n(x)||,||z^n_i||-1\}.$
Otherwise, if $x\in D_n\setminus M_{n,i}$ we proceed by induction in $c(n)-i\in\N$. For the first step, $i=c(n)-1$ and $x\in D_n\setminus M_{n,c(n)-1}$ so that $x=c^n_{c(n)}$. Then,
$$\begin{aligned}||r_n\circ \phi_{n,c(n)-1}(x)||=&||r_n\circ T_{n,c(n)}(c^n_{c(n)})||=||r_n(c^n_{c(n)})||-1=||z_{c(n)}^n||-1\\\ge&||z_{c(n)-1}^n||-1\ge\min\{||r_n(x)||,||z_{c(n)-1}^n||-1\}.\end{aligned}$$
To prove the inductive step we take $i\in\{0,\dots,c(n)-2\}$ and assume by induction hypothesis that
\begin{equation}\label{step22}||r_n\circ \phi_{n,i+1}(x)||\ge \min\{ ||r_n(x)||,||z^n_{i+1}||-1 \}.\end{equation}
Now, \eqref{step2} follows from \eqref{ineqind} and \eqref{step22} since
$$\begin{aligned}||r_n\circ \phi_{n,i}(x)||=&||r_n\circ T_{n,i+1}\circ \phi_{n,i+1}(x)||\ge\min\{||r_n\circ \phi_{n,i+1}(x)||,||z^n_i||-1\}\\\ge& \min\{||r_n(x)||,||z^n_{i+1}||-1 ,||z^n_{i}||-1 \}=\min\{||r_n(x)||,||z^n_{i}||-1 \}.\end{aligned}$$
Now that we know \eqref{step2} holds true, we may compute the Lipschitz norm of $\phi_{n,i}$. Let us take distinct points $x,y\in D_{n}$  (meaning that $r_n(x)\neq r_n(y)$), assuming without loss of generality that $x\notin M_{n,i}$ (otherwise the statement is trivially satisfied). Since $x\notin M_{n,i}$ we know that $||r_n(x)||\ge||z_i^n||$. By \eqref{step2} it follows that
$$||r_n\circ \phi_{n,i}(y)||-||r_n\circ \phi_{n,i}(x)||\le ||z^n_i||-\min\{||r_n(x)||,||z_i^n||-1\}=1\le||x-y||.$$
Also, taking again into account that $||z_i^n||\le||r_n(x)||$,
$$||r_n\circ \phi_{n,i}(x)||-||r_n\circ \phi_{n,i}(y)||\le ||z^n_i||-\min\{||r_n(y)||,||z^n_i||-1\}\le||x-y||.$$
Therefore, using Lemma \ref{defstep4} we obtain that
$$\begin{aligned}||r_n\circ \phi_{n,i}(x)-r_n\circ \phi_{n,i}(y)||=&||T_{||r_n\circ \phi_{n,i}(x)||}(r_n(x))-T_{||r_n\circ \phi_{n,i}(y)||}(r_n(y))||\\\le&||T_{||r_n\circ \phi_{n,i}(x)||}(r_n(x))-T_{||r_n\circ \phi_{n,i}(y)||}(r_n(x))||\\&+||T_{||r_n\circ \phi_{n,i}(y)||}(r_n(x))-T_{||r_n\circ \phi_{n,i}(y)||}(r_n(y))||\\\le& \big|\,||r_n\circ \phi_{n,i}(x)||-||r_n\circ \phi_{n,i}(y)||\,\big|\\&+||r_n(x)-r_n(y)||\\\le&2||x-y||.\end{aligned}$$
Hence,
$$\begin{aligned}||\phi_{n,i}(x)-\phi_{n,i}(y)||\le& ||\big( \restr{r_n}{D_n} \big)^{-1}\big(r_n\circ \phi_{n,i}(x)\big)-\big( \restr{r_n}{D_n} \big)^{-1}\big(r_n\circ \phi_{n,i}(y)\big)||\\\le&(\lambda^2+2\lambda+2)||r_n\circ \phi_{n,i}(x)-r_n\circ \phi_{n,i}(y)||\\\le&(\lambda^2+2\lambda+2)2||x-y||.\end{aligned}$$
\end{proof}
\end{prop}

We are now done with the forth step,  proving Proposition \ref{step4} with $L=2\lambda^2+4\lambda+4$.

\subsection{Step 5}

In this step we are finally constructing a retractional basis for $M$. It is just a matter of amalgamating the retractions of Proposition \ref{step3} and Proposition \ref{step4} together.

\begin{theorem}
$M$ has a retractional basis. 
\begin{proof}
We list the elements of $M_1$ as $(x^1_1,\dots,x^1_{i(1)})$ where $i(1)=\# M_1$ and $x_1^1=0\in M_1$. Let us define the mapping $I:M\to \N$ as
$$I(x)=\begin{cases}i\;\;&\text{if }x=x^1_i\in M_1,\\i+\#M_n\;\;&\text{if }x=c^n_i\in C_n=D_n\setminus M_n\text{ for }n\in\N,\\i+\#D_n\;\;&\text{if }x=x^{n}_i\in M_{n}\setminus D_{n-1}\text{ for }n\ge2.\end{cases}$$
$I$ is clearly well defined and bijective. It is straightforward to see that
\begin{itemize}
\item $I(x)\in\{1,\dots,i(1)\}$ if, and only if, $x=x^1_{I(x)}\in M_1.$
\item $I(x)\in\{\#M_n+1,\dots,\#D_n\}$ if, and only if, $x=c^n_{I(x)-\#M_n}\in D_{n}\setminus M_n=C_n.$
\item $I(x)\in\{\#D_{n-1}+1,\dots,\#M_n\}$ if, and only if, $x=x^n_{I(x)-\#D_{n-1}}\in M_n\setminus D_{n-1}.$
\end{itemize}
We denote $M^i=(I^{-1}(l))_{l\le i}$ for every $i\in\N$. For each $i\in\N$ the retraction 
of the retractional basis onto $M^i$  $\varphi_i:M\to M^i$ is given by
$$\varphi_i(x)=\begin{cases}x\;\;&\text{if }x\in M^i,\\\phi_{n,i-\#M_n}\circ \Psi_{n+1}\circ \phi_{n+1}(x)\;\;&\text{if } I^{-1}(i)\in D_n\setminus M_n\text{ for }n\in\N,\\\Psi_{n,i-\#D_{n-1}}\circ\phi_n(x)\;\;&\text{if } I^{-1}(i)\in M_n\setminus D_{n-1}\text{ for }n\ge2,\\0\;\;&\text{if }I^{-1}(i)\in M_1,\; x\notin M^i. \end{cases}$$
An easy case by case computation shows that $\varphi_i$ is a well defined retraction onto $M^i$ satisfying $\varphi_{i_1}\circ\varphi_{i_2}=\varphi_{\min(i_1,i_2)}$ for every $i_1,i_2\in\N$. In fact, it is enought to use Proposition \ref{step3} and Proposition \ref{step4} together with the fact that for every $n_1,n_2\in\N$ with $n_2\ge2$ and $n_2\ge n_1$,
$$\phi_{n_1}=\phi_{n_1}\circ \Psi_{n_2,i_2}\circ \phi_{n_2+1}.$$

It only remains to prove that $\varphi_i$ is Lipschitz. Let us take distinct $x,y\in M$. The case $x,y\in M^i$ is trivial since $\varphi_i(x)-\varphi_i(y)=x-y$. Hence, we may assume that $x\notin M^i$. If $I^{-1}(i)\in D_n\setminus M_n=C_n$ then $\varphi_i(x)=\phi_{n,i-\#M_n}\circ \Psi_{n+1}\circ \phi_{n+1}(x)$ and $\varphi_i(y)=\phi_{n,i-\#M_n}\circ \Psi_{n+1}\circ \phi_{n+1}(y)$. Therefore,
$$\begin{aligned}||\varphi_i(x)-\varphi_i(y)||\le&||\phi_{n,i-\#M_n}||_{\text{Lip}}||\Psi_{n+1}||_{\text{Lip}}||\phi_{n+1}||_{\text{Lip}}||x-y||\\\le& (2\lambda^2+4\lambda+4)(\lambda^2+2\lambda+2)(3\lambda+2)||x-y||.\end{aligned}$$
If $I^{-1}(i)\in M_n\setminus D_{n-1}$ then $\phi_i(x)=\Psi_{n,i-\#D_{n-1}}\circ\phi_n(x)$ and $\phi_i(y)=\Psi_{n,i-\#D_{n-1}}\circ\phi_n(y)$. Therefore,
$$\begin{aligned}||\varphi_i(x)-\varphi_i(y)||\le&||\Psi_{n,i-\#D_{n-1}}||_{\text{Lip}}||\phi_{n}||_{\text{Lip}}||x-y||\\\le&(3\lambda+2)^2(3\lambda^2+6\lambda+11)||x-y||.\end{aligned}$$
Finally, we focus on the case when $I^{-1}(i)\in M_1$. If $y\in M^i$ then $y\in M^1=f\circ i_n\circ f(\lambda B_{\ell_\infty(\Gamma_1)})$ so $||y||\le\lambda^2$. Hence, $||\varphi_i(x)-\varphi_i(y)||=||y||\le\lambda^2\le\lambda^2||x-y||$. Otherwise, if $y\notin M^i$ the Lipschitz inequality is trivially satisfied since $\phi_i(x)=\phi_i(y)=0$.


\end{proof}
\end{theorem}

In general $M$ is not a net of $X$ because it may not be contained in $X$. 
Replacing it with a net in $X$ with the same essential properties
 is the aim of the last step.

\subsection{Step 6}

We first prove the next simple lemma.

\begin{lemma}\label{densityM}
For every $x\in X$ there is $y\in M$ such that $||x-y||\le 2\lambda+2$.
\begin{proof}
Clearly, for any $x\in X$ there exists $n\in\N$ such that $||x-i_n\circ r_n(x)||\le1$ because $(i_n\circ r_n)_{n\in\N}$ are the natural projections of a Schauder Decomposition of $X$ (Remark 2.9 of \cite{AGM16}). Therefore, $\widetilde x:=i_n\circ r_n\circ f(x)\in\widetilde M_n$ satisfy that $||x-\widetilde x||\le||x-i_n\circ r_n(x)||+||i_n\circ r_n(x)-i_n\circ r_n\circ f(x)||\le 1+\lambda$. We let $y=\big(\restr{r_n}{M_n}\big)^{-1}(r_n(\widetilde x))$ so that $i_n\circ r_n(y)=i_n\circ r_n(\widetilde x)=\widetilde x$. By Lemma \ref{dense} we have that
$$||\widetilde x-y||=||i_n\circ r_n(y)-y||\le \lambda+1.$$
Finally,
$$||x-y||\le ||x-\widetilde x||+||\widetilde x-y||\le2\lambda+2,$$
and we are done.
\end{proof}
\end{lemma}

\begin{prop}\label{LIPnet}
There is a net $N$ in $X$ which is Lipschitz equivalent to $M$.
\begin{proof}
For every $m\in M$ there is a unique $n(m)\in\N$ such that $m\in M_{n(m)}\setminus M_{n(m)-1}$ (where $M_0=\emptyset$). We define  the map $\rho:M\to X$ as $\rho(m)=i_{n(m)}\circ r_{n(m)}(m)$ for every $m\in M$. If $m\in M$, it follows from the definition of $M_n\setminus M_{n-1}$ that $m=f\circ \rho(m)$ so $||m-\rho(m)||\le1$. The problem here is that $\rho$ may not be injective and $\rho(M)$ may not be uniformly discrete so we are going to force it. First, we order $\rho(M)$ (which is countable) as $(m_n)_{n\in\N}$ and produce a sequence of points $(\widetilde m_n)\subset \rho(M)$ together with a sequence of disjoint non-empty subsets $N_n\subset \rho(M)$ inductively. Let us take $a>1$, we set $\widetilde m_1=m_1$ and $N_1=\{m\in\rho(M)\;:\;||m-\widetilde m_1||\le a\}$ and if $\widetilde m_k$ and $N_k$ are known for every $k<n$ we set $\widetilde m_n=\min \rho(M)\setminus \big( N_1\cup\cdots\cup N_{n-1} \big)$ and $N_n=\{m\in\rho(M)\setminus\big(N_1\cup\cdots\cup N_{n-1}\big)\;:\;||m-\widetilde m_n||\le a\}$. Obviously, $\widetilde N:=(\widetilde m_n)$ is an $(a,a)$-net of $\rho(M)$. Hence, it is an $(a,b)$-net of X with $b=a+2\lambda+3$, since  by Lemma \ref{densityM}, $\rho(M)$ is $(2\lambda+3)$-dense in $X$. Also, $\rho(M)=\bigcup\limits_{k\in\N}N_k$ where $N_i\cap N_j=\emptyset$ for every $i,j\in\N$ with $i\neq j$.

Now, since $X$ is infinite dimensional, there is a sequence of points $(x_n)\subset \frac{3}{4}B_X$ such that $||x_m-x_n||\ge 1/2$ for every $m,n\in\N$ with $m\neq n$. To simplify the next definition, we are going to  order the points of $\rho^{-1}(N_k)=(m_k^n)_n$ for every $k\in\N$ (this is possible because $M$ is countable). Finally, since $M=\bigcup\limits_{k\in\N}\rho^{-1}(N_k)$, we define the injection $\mu:M\to X$ as
$$\mu(m_k^n)=\widetilde{m}_k+\frac{a}{2} x_n.$$
We claim that $N=\mu(M)$ is a net of $X$ Lipschitz equivalent to $M$. Clearly, $N$ is $(b+3a/8)$-dense in $X$ since $\widetilde N$ is $b$-dense in $X$. It is also easy to see that
$$||\mu(m_{k_1}^{n_1})-\mu(m_{k_2}^{n_2})||\ge ||\widetilde m_{k_1}-\widetilde m_{k_2}||-\frac{a}{2}\big(||x_{n_1}||+||x_{n_2}||\big)\ge a-\frac{a}{2}\bigg( \frac{3}{4}+\frac{3}{4} \bigg)= a/4$$
whenever $k_1\neq k_2$. If $k_1=k_2$ then obviously $||\mu(m_{k_1}^{n_1})-\mu(m_{k_2}^{n_2})||=a/2||x_{n_1}-x_{n_2}||\ge a/4$. This proves that $N$ is a $(a/4,b+3a/8)$-net in $X$. It only remains to prove that $\mu$ and $\mu^{-1}$ (seeing $\mu$ restricted to $N$) are Lipschitz. This follows immediately from the fact that for every $m_k^n\in M$,
$$\begin{aligned}||\mu(m_k^n)-m_{k}^n||=&||\widetilde{m}_k-m_{k}^n+x_na/2||\\\le&||\widetilde m_k-\rho(m_k^n)||+||\rho(m_k^n)-m_k^n||+||x_n a/2||\\\le& a+1+3a/8.\end{aligned}$$
\end{proof}
\end{prop}

\begin{theorem}[Main Theorem]
A net in any  separable infinite dimensional $\mathcal{L}_\infty$-space has a retractional basis.
\end{theorem}

\begin{corollary}
If $N$ is a  net in a separable infinite dimensional $\mathcal{L}_\infty$-space 
then $\mathcal{F}(N)$ has a Schauder basis.
\end{corollary}

\bigskip

\printbibliography
\end{document}